\RequirePackage{ifpdf}
\ifpdf 
\documentclass[pdftex]{sigma}
\else
\documentclass{sigma}
\fi

\usepackage{mathrsfs}

\newcommand{\p}{\partial}
\begin{document}

\allowdisplaybreaks

\renewcommand{\PaperNumber}{044}

\FirstPageHeading

\ShortArticleName{Generalized Hasimoto Transform of
One-Dimensional Dispersive Flows}

\ArticleName{Generalized Hasimoto Transform of One-Dimensional\\
Dispersive Flows into Compact Riemann Surfaces}

\Author{Eiji ONODERA}

\AuthorNameForHeading{E. Onodera}

\Address{Mathematical Institute, Tohoku University, Sendai
980-8578, Japan}
\Email{\href{mailto:sa3m09@math.tohoku.ac.jp}{sa3m09@math.tohoku.ac.jp}}

\ArticleDates{Received December 18, 2007, in f\/inal form May 14,
2008; Published online May 20, 2008}

\Abstract{We study the structure of dif\/ferential equations of
one-dimensional dispersive f\/lows into compact Riemann surfaces.
These equations geometrically generalize two-sphere valued systems
modeling the motion of vortex f\/ilament. We def\/ine a
generalized Hasimoto transform by constructing a good moving
frame, and reduce the equation with values in the induced bundle
to a complex valued equation which is easy to handle. We also
discuss the relationship between our reduction and the theory of
linear dispersive partial dif\/ferential equations.}


\Classification{35Q55; 35Q35; 53Z05}

\section{Introduction}
\label{section:introduction} We study the structure of dispersive
partial dif\/ferential equations of the form
\begin{gather}
  u_{t}
 =
  a\,\nabla_x^2u_x
  + J_{u}\nabla_xu_x
  +b\,g_{u}(u_x, u_x)u_x,
\label{equation:pde3}
\\
   u_t
  =
   -
   a\, J_u\nabla_x^3u_x
   +
   \{ 1+b\, g_u(u_x,u_x) \}J_u\nabla_xu_x
   +
   c\, g_u(\nabla_xu_x,u_x)J_uu_x,
\label{equation:pde4}
\end{gather}
where $N$ is a Riemann surface with an almost complex structure
$J$ and a hermitian metric $g$, $a, b, c \in \mathbb{R}$ are
constants, $u=u(t,x)$ is an $N$-valued unknown function of
$(t,x)\in \mathbb{R}\times \mathbb{R}$, $u_t(t,x) =du_{(t,x)} (
\left(\p/\p{t}\right)_{(t,x)})$, $u_x(t,x) =du_{(t,x)}(
\left(\p/\p{x}\right)_{(t,x)})$,
$du_{(t,x)}:T_{(t,x)}\mathbb{R}^2\to T_{u(t,x)}N$ is the
dif\/ferential of the mapping $u$ at $(t,x)$, $\nabla_x$ is the
covariant derivative with respect to $x$ along the mapping $u$,
and $J_u$ and $g_u$ denote the almost complex structure and the
metric at $u \in N$ respectively. Recall that an almost complex
structure $J$ on $N$ is a bundle automorphism of the tangent
bundle $TN$ satisfying $J^2=-I$. Note that~\eqref{equation:pde3}
and~\eqref{equation:pde4} are equations for vector f\/ields along
the mapping $u$.

In mathematics, a non-Kowalewskian-type evolution equation is said
to be ``dispersive'' if its initial value problem is (expected to
be) well-posed in both directions of the past and the future. See
\cite[Lecture VII]{Mizohata}. Loosely speaking, in classical
mechanics, a linear evolution equation with constant
coef\/f\/icients is called ``dispersive'' if the propagation speed
of its plane wave solutions changes in proportion to the size of
frequency. Typical examples of dispersive partial dif\/ferential
equations are the Schr\"odinger evolution equation of free
particle $v_t+iv_{xx}=0$, $i=\sqrt{-1}$ and the Korteweg--de Vries
equation $v_t+v_{xxx}-6vv_x=0$. We call the solution of
\eqref{equation:pde3} or \eqref{equation:pde4} a dispersive
f\/low. In particular, when $a=b=0$, the solution of
\eqref{equation:pde3} is called a one-dimensional Schr\"odinger
map.

Dispersive f\/lows arise in classical mechanics: the motion of
vortex f\/ilament, the Heisenberg ferromagnetic spin chain and
etc. See \cite{DR,Fukumoto, FM,Hama,Hasimoto,SSB} and references
therein. Solutions to such equations in classical mechanics are
valued in two-dimensional unit sphere. The
equation~\eqref{equation:pde3} is the geometric generalization of
equations proposed by Fukumoto and Miyazaki in~\cite{FM},
and~\eqref{equation:pde4} is the geometric generalization of
equations formulated by Fukumoto in~\cite{Fukumoto}. See
Section~\ref{section:modeling} for details.

The replacement of the sphere as a target by a more general
Riemannian manifold is one of the ways to study these physical
models with some geometric structures like~\eqref{equation:pde3}
or~\eqref{equation:pde4}. In fact, there has been many geometric
analyses on the relationship between
 the good structure of~\eqref{equation:pde3} or~\eqref{equation:pde4}
and the geometry of $N$, where $N$ is a compact Riemann surface,
or more generally, a~compact K\"ahler manifold. See, e.g.,
\cite{CSU,Koiso,Koiso2,Onodera} and~\cite{PWW}.

Especially, the equation of Schr\"odinger map geometrically
generalizes an equation proposed by Da Rios in \cite{DR}. In
\cite{CSU} Chang, Shatah and Uhlenbeck studied the structure of
\eqref{equation:pde3} with $a=b=0$. They constructed a moving
frame of the induced bundle $u^{-1}TN$, and reduced the equation
to a simple form of complex valued semilinear partial
dif\/ferential equation. More precisely, if there exists
$u^\ast\in{N}$ such that
\begin{gather*}
\lim_{x\rightarrow-\infty}u(t,x)=u^\ast
\end{gather*}
for any $t \in \mathbb{R}$, \eqref{equation:pde3} with $a=b=0$ is
reduced to
\begin{gather}
q_t=iq_{xx}+\frac{i}{2}\kappa(u)\left|q \right|^2q
    -\frac{i}{2}
    \left[
    \int_{-\infty}^x
    \left|q \right|^2
    \left( \kappa(u)\right)_x
    dx'
     \right]q,
\label{equation:hasimoto2}
\end{gather}
where $q(t, x)$ is a complex-valued function of $(t,x)$ and
$\kappa(u)$ denotes the Gaussian curvature at $u\in N$. Their idea
of this reduction came from the work of Hasimoto in
\cite{Hasimoto}. Their idea is also essentially related with the
basic method for one-dimensional linear Schr\"odinger-type
evolution equations. See Section~\ref{section:remarks} for
details.

The purpose of this paper is to reduce the equations of dispersive
f\/lows \eqref{equation:pde3} and \eqref{equation:pde4} to complex
valued equations by constructing a good moving frame. To state our
results, we introduce Sobolev spaces $H^\infty(\mathbb{R})$ and
$H^\infty(\mathbb{R};TN)$ def\/ined by
\begin{gather*}
  H^\infty(\mathbb{R})
 =
  \left\{
  q{\in}C^\infty(\mathbb{R};\mathbb{C})
  \ \bigg\vert \
  \int_{\mathbb{R}}
  \big\lvert \p^k_xq (x)\big\rvert^2
  dx<\infty,
  \ k=0,1,2,\dots
  \right\},
\\
  H^\infty(\mathbb{R};TN)
 =
  \left\{
  u{\in}C^\infty(\mathbb{R};N)
  \ \bigg\vert \
  \int_{\mathbb{R}}
  g(\nabla_x^ku_x,\nabla_x^ku_x)
  dx<\infty,
  \ k=1,2,\dots
  \right\}.
\end{gather*}
For a function space $X$, $C(\mathbb{R};X)$ denotes the set of all
$X$-valued continuous functions on $\mathbb{R}$. Our results are
the following.

\begin{theorem}
\label{theorem:hasimoto3} Assume that the equation
\eqref{equation:pde3} has a solution $u\in
C(\mathbb{R};H^\infty(\mathbb{R};TN))$ and that $\lim\limits_{x\to
-\infty}u(t,x)= u^{*}\in N$ for any $t\in\mathbb{R}$. Then there
exists a complex-valued function $q\in  C(\mathbb{R};
H^{\infty}(\mathbb{R}))$ solving
\begin{gather}
 q_t -aq_{xxx}-iq_{xx}
  =
 \left(
  \frac{a}{2}\kappa(u)+2b\right)
  \lvert{q}\rvert^{2}
  q_{x}
  -
  \left(
  \frac{a}{2}\kappa(u)-b
  \right)
  q^{2}\bar{q}_{x}\nonumber
\\
\qquad{}
 +
 ia
  \left[
  \int_{-\infty}^{x}
  (\kappa(u))_{x}
  \operatorname{Im}\,(q\bar{q}_{x})dx'
  \right]
  q
  -
 \frac{i}{2}
  \left[
  \int_{-\infty}^{x}
  (\kappa(u))_{x}
  \lvert{q}\rvert^{2}dx'
  \right]
  q
  +
  \frac{i}{2}
  \kappa(u)
  \lvert{q}\rvert^{2}q.
\label{equation:t3rd}
\end{gather}
\end{theorem}

\begin{theorem}
\label{theorem:hasimoto4} Assume that the equation
\eqref{equation:pde4} has a solution $u\in
C(\mathbb{R};H^\infty(\mathbb{R};TN))$ and that $\lim\limits_{x\to
-\infty}u(t,x)= u^{*}\in N$ for any $t\in\mathbb{R}$. Then there
exists a complex-valued function $q\in  C(\mathbb{R};
H^{\infty}(\mathbb{R}))$ solving
\begin{gather}
 q_t +ia q_{xxxx}-iq_{xx}
  =
 i
  \left(
  b
  +
  \frac{c}{2}
  -
  \frac{a}{2}\kappa(u)
  \right)
  \lvert{q}\rvert^2
  q_{xx}
  +
  i
  \left(
  \frac{c}{2}
  -
  \frac{a}{2}\kappa(u)
  \right)
  q^2\bar{q}_{xx}
  +
  i
  \left(
  b
  +
  \frac{c}{2}
  \right)
  \bar{q}q_{x}^2
\nonumber
\\
\qquad{} +i
  \left(
  b
  +
  \frac{3c}{2}
  +
  \frac{a}{2}\kappa(u)
  \right)
  q\lvert{q_x}\rvert^2
  +
  \frac{i}{2}
  \kappa(u)
  \lvert{q}\rvert^2q
  +i
  \frac{b+c}{4}
  \kappa(u)
  \lvert{q}\rvert^4q
\nonumber
\\
\qquad{} -
  i
  \left[
  \int_{-\infty}^x
  \Bigl(\kappa(u)\Bigr)_x
  \left\{
  \frac{1}{2}
  \lvert{q}\rvert^2
  +
   \frac{b+c}{4}
  \lvert{q}\rvert^4
  +
  \frac{a}{2}
  \lvert{q_x}\rvert^2
  -
  \frac{a}{2}
  q_{xx}\bar{q}
  -
   \frac{a}{2}
  \bar{q}_{xx}q
  \right\}
  dx'
  \right]
  q.
\label{equation:t4th}
\end{gather}
\end{theorem}

The meaning of our reductions is to see the essential structure of
the equations \eqref{equation:pde3} and \eqref{equation:pde4} from
a view point of the theory of linear partial dif\/ferential
operators. Using our reductions, one can easily understand whether
the classical method of solving initial value problems works or
not. See Section~\ref{section:remarks} for details. We remark that
our reductions cannot be used as the tool of solving the initial
value problem for the equations~\eqref{equation:pde3}
and~\eqref{equation:pde4}. Firstly, we impose that solutions have
a~f\/ixed base point $u^\ast$ at $x=-\infty$, and such solutions
are not physically interesting. Secondly, our
reductions~\eqref{equation:t3rd} and~\eqref{equation:t4th} contain
functions depending on solutions $u$ to the original
equations~\eqref{equation:pde3} and~\eqref{equation:pde4}, and we
need to express $u$ by $q$ to solve \eqref{equation:t3rd} and
\eqref{equation:t4th}. Our reductions can be used as the tool of
solving \eqref{equation:pde3} and \eqref{equation:pde4} only in
the case that the maps are near constant and the target is the
two-sphere or a hyperbolic surface.

The plan of this article is as follows. In
Section~\ref{section:modeling} we derive \eqref{equation:pde3} and
\eqref{equation:pde4} from classical mechanical models. In
Section~\ref{section:proofs} we prove
Theorems~\ref{theorem:hasimoto3} and \ref{theorem:hasimoto4}.
Finally, in Section~\ref{section:remarks} we discuss our reduction
from a point of view of the theory of linear dispersive partial
dif\/ferential equations. Throughout this paper, as has already
been used in the above equation \eqref{equation:hasimoto2},
\eqref{equation:t3rd} and \eqref{equation:t4th}, the partial
dif\/ferentiation for $\mathbb{R}$ or $\mathbb{C}$-valued
functions is written by $\p$ or the script, e.g., $\p_xf$, $f_x$.

\section{Geometric generalization of physical models}
\label{section:modeling} Consider the motion of a very thin
isolated vortex f\/ilament in three-dimensional incompressible
unbounded perfect f\/luid. We denote the position of the vortex
f\/ilament by $\vec{\bf{X}}(t,x)$, where $t$ is a~time, and $x$ is
the arc length. Let $(\kappa,\tau)$ be the curvature and the
torsion and let $(\vec{\bf{X}}_x, \vec{\bf{n}}, \vec{\bf{b}})$ be
the Frenet--Serret frame along the vortex f\/ilament
$\vec{\bf{X}}$. In \cite{DR} Da Rios formulated the following
model equation
\begin{gather}
\vec{\bf{X}}_t = \kappa\vec{\bf{b}}. \label{equation:vortex1}
\end{gather}
by using the so called localized induction approximation. See
also~\cite{Hama}. In~\cite{Hasimoto} Hasimoto pointed out that
\eqref{equation:vortex1} can be formally transformed into the
cubic nonlinear Schr\"odinger equation
\begin{gather}
\psi = i\psi_{xx}+ \frac{i}{2}\left[ |\psi|^2+A(t) \right]\psi
\label{equation:vortex4}
\end{gather}
by using
\begin{gather}
\psi(t,x) = \kappa(t,x) \exp \left( i\int_{0}^x \tau(t,x')dx'
\right), \label{equation:hasimoto}
\end{gather}
where $A(t)$ is a real-valued function of $t\in \mathbb{R}$.
Furthermore, if we set
\begin{gather*}
\phi(t,x) = \psi(t,x) \exp\left(
-\frac{i}{2}\int_{0}^tA(t')dt'\right),
\end{gather*}
\eqref{equation:vortex4} becomes
\begin{gather*}
\phi_t=i\phi_{xx}+\frac{i}{2}|\phi|^2\phi.
\end{gather*}
The correspondence $X{\mapsto}\psi$ or $X{\mapsto}\phi$ is said to
be the Hasimoto transform.

Recently, other model equations of vortex f\/ilament were
proposed. In \cite{FM} Fukumoto and Miyazaki proposed the
following
\begin{gather}
\vec{\bf{X}}_t = \kappa\vec{\bf{b}}+ a \left[
\frac{1}{2}\kappa^2\vec{\bf{X}}_x + \kappa_x\vec{\bf{n}} +
\kappa\tau\vec{\bf{b}} \right], \label{equation:vortex5}
\end{gather}
where the real constant $a\in \mathbb{R}$ is the ef\/fect of axial
f\/low f\/lux through the vortex tube. In \cite{Fukumoto},
considering the inf\/luence of elliptical deformation of the core
due to the self-induced strain, Fukumoto also proposed the
following
\begin{gather}
\vec{\bf{X}}_t = \kappa\vec{\bf{b}}+ C_1\left[
\kappa^2\tau\vec{\bf{X}}_x + \left(2\kappa_x\tau+\kappa\tau_x
\right)\vec{\bf{n}} + \left( \kappa\tau^2-\kappa_{xx}
\right)\vec{\bf{b}} \right] +C_b\kappa^3\vec{\bf{b}},
\label{equation:vortex6}
\end{gather}
where $C_1$, $C_b$ are real constants.

\eqref{equation:vortex5} can be transformed by the Hasimoto
transform into the following equation of the form
\begin{gather}
\psi_t= i\psi_{xx}+\frac{i}{2}|\psi|^2\psi +a \left\{
\psi_{xxx}+\frac{3}{2}|\psi|^2\psi_x \right\}.
\label{equation:hirota}
\end{gather}
This equation \eqref{equation:hirota} is called a Hirota equation,
and is known as the completely integrable equation in classical
mechanics. Similarly, \eqref{equation:vortex6} can be transformed
to the following equation of the form
\begin{gather}
  \psi_t
 =
  i\psi_{xx}+\frac{i}{2}\lvert\psi\rvert^2\psi
 -
  iC_1
  \left\{
  \psi_{xxxx}
  +
  \frac{3}{2}
  \bigl(
  \lvert\psi\rvert^2\psi_{xx}
  +
  \psi_x^2\bar{\psi}
  \bigr)
  +
  \left(
  \frac{3}{8}\lvert\psi\rvert^4
  +
  \frac{1}{2}
  \bigl(\lvert\psi\rvert^2\bigr)_{xx}
  \right)\psi
  \right\}
\nonumber
\\
\phantom{\psi_t=}{} +
  i
  \left(C_b+\frac{C_1}{2}\right)
  \left\{
  \bigl(\lvert\psi\rvert^2\psi\bigr)_{xx}
  +
  \frac{3}{4}\lvert\psi\rvert^4\psi
  \right\}.
\label{equation:4th}
\end{gather}

We show that the equations \eqref{equation:vortex5} and
\eqref{equation:vortex6} can be written in the form of equations
\eqref{equation:pde3} and \eqref{equation:pde4}, respectively. For
$\vec{u}=(u_1,u_2,u_3)\in\mathbb{R}^3$ and
$\vec{v}=(v_1,v_2,v_3)\in\mathbb{R}^3$, set
\begin{gather*}
(\vec{u},\vec{v})=u_1v_1+u_2v_2+u_3v_3, \qquad
\lvert\vec{u}\rvert=\sqrt{(\vec{u}, \vec{u})   },
\\
\vec{u}\times\vec{v} =
(u_2v_3-u_3v_2,u_3v_1-u_1v_3,u_1v_2-u_2v_1),
\end{gather*}
and denote $\mathbb{S}^2 = \{\vec{u}\in\mathbb{R}^3 \ \vert \
\lvert\vec{u}\rvert=1\}$. If the curvature of $\vec{\bf{X}}$ never
vanishes, that is $\vec{\bf{X}}_{xx}\neq\vec{0}$, the right hand
side of \eqref{equation:vortex5} and \eqref{equation:vortex6} make
sense and are written as
\begin{gather}
\vec{\bf{X}}_t = \vec{\bf{X}}_x \times \vec{\bf{X}}_{xx} + a
\left[ \vec{\bf{X}}_{xxx} + \frac{3}{2}
\lvert\vec{\bf{X}}_{xx}\rvert^2 \vec{\bf{X}}_x \right],
\label{equation:vortex7}
\\
  \vec{\bf{X}}_t
 =
  \vec{\bf{X}}_x
  \times
  \vec{\bf{X}}_{xx}
  -
  C_1
  \vec{\bf{X}}_x
  \times
  \vec{\bf{X}}_{xxxx}
 +
  C_1
  \vec{\bf{X}}_{xx}
  \times
  \vec{\bf{X}}_{xxx}
  +
  (C_b-2C_1)|\vec{\bf{X}}_{xx} |^2
  \vec{\bf{X}}_x
  \times
  \vec{\bf{X}}_{xx},
\label{equation:vortex8}
\end{gather}
respectively. Dif\/ferentiating the vortex f\/ilament equation
\eqref{equation:vortex7} and \eqref{equation:vortex8} with respect
to $x$, we obtain the equation for the velocity vector
$\vec{\bf{u}}=\vec{\bf{X}}_x$ of the form
\begin{gather}
\vec{\bf{u}}_t = \vec{\bf{u}} \times \vec{\bf{u}}_{xx} + a \left\{
\vec{\bf{u}}_{xxx} + 3(\vec{\bf{u}}_{xx},
\vec{\bf{u}}_x)\vec{\bf{u}} + \frac{3}{2} \left|\vec{\bf{u}}_x
\right|^2 \vec{\bf{u}}_x \right\}, \label{equation:vortex9}
\\
\vec{\bf{u}}_t = \vec{\bf{u}} \times \vec{\bf{u}}_{xx} - C_1
\vec{\bf{u}} \times \vec{\bf{u}}_{xxxx} + (C_b-2C_1) \big(
\left|\vec{\bf{u}}_x \right|^2\vec{\bf{u}} \times  \vec{\bf{u}}_x
\big)_x, \label{equation:vortex10}
\end{gather}
respectively. Since $x$ is the arc length,
$\lvert\vec{\bf{u}}\rvert=1$, which implies that $\vec{\bf{u}}$
lies on $\mathbb{S}^2$.

Here assume that $u(t,x)\in \mathbb{S}^2$ for any $(t,x)$. Let $V$
be  any  vector f\/ield along $u$, namely, $V(t,x)\in
T_u\mathbb{S}^2$ for any $(t,x)$. Regarding $\mathbb{S}^2$ as the
two-sphere with the standard metric induced from $\mathbb{R}^3$,
we see from the def\/inition of the covariant derivative on
$\mathbb{S}^2$ that
\begin{gather}
\nabla_xV=V_x-(V_x, u)u=V_x+(V, u_x)u. \nonumber
\end{gather}
Thus we have
\begin{gather*}
\nabla_xu_x=u_{xx}+\left|u_x \right|^2u,
\\
\nabla_x^2u_x=u_{xxx}+3\left(u_{xx}, u_x\right)u
                      +\left|u_x \right|^2u_x,
\\
\nabla_x^3u_x=u_{xxxx}+4\left(u_{xxx}, u_x\right)u
                      +3\left|u_{xx} \right|^2u
                      +5\left(u_{xx}, u_x\right)u_x
                      +\left|u_{x} \right|^2u_{xx}
                      +\left|u_{x} \right|^4u.
\end{gather*}
Also, the operation $u\times V$, rotating $V$ by $\pi/2$ degrees
for any $V{\in}T_{u}\mathbb{S}^2$, acts as a complex number in the
tangent space. Thus $J_{u}V=u\times V$ for any vector f\/ield $V$
along $u$, which implies
\begin{gather}
u\times u_{xx}=u\times \nabla_xu_x=J_u\nabla_xu_x. \nonumber
\end{gather}
Noting these relations, we see that \eqref{equation:vortex9} is
written as \eqref{equation:pde3} with $b=a/2$ and
\eqref{equation:vortex10} is written as \eqref{equation:pde4} with
$a=C_1$, $b=C_b-C_1$, $c=2C_b+C_1$ respectively.

\section{Proofs}
\label{section:proofs} In this section, we prove
Theorem~\ref{theorem:hasimoto3} and \ref{theorem:hasimoto4}
proceeding as in~\cite{CSU}.

\begin{proof}[Proof of Theorem~\ref{theorem:hasimoto3}]
Assume that there exists a point $u^{*}\in N$ and that equation
\eqref{equation:pde3} has a~smooth solution on $\mathbb{R} \times
\mathbb{R}$ such that $u(t,x)\to u^{*}\in N$ as $x\to -\infty$.
Let $\left\{e, Je \right\}$ denote the orthonormal frame for
$u^{-1}TN$ constructed in the following manner:

Fix a unit vector $e_0\in T_{u^{*}}N$, namely, $g_{u^{*}}(e_0,
e_0)=1$. And for any $t\in \mathbb{R}$, let $e(t,x)\in
T_{u(t,x)}N$ be the parallel translation of $e_0$ along the curve
$u(t,\cdot )$, that is,
\begin{gather}
\nabla_xe(t,x)=0, \label{equation:frame1}
\\
\lim_{x\to -\infty}e(t,x)=e_0. \label{equation:bound}
\end{gather}
In local coordinate $(U, \phi =(u^1, u^2))$ around $u(t,x)$, where
the metric is given by $g=\lambda(z,\bar{z})dzd\bar{z}$ with a
complex coordinate $z=u^1+iu^2$, set
\begin{gather*}
e(t,x)=e^1(t,x)\left(\frac{\p}{\p u^1}\right)_u
+e^2(t,x)\left(\frac{\p}{\p u^2}\right)_u, \qquad \text{and}
\qquad \zeta =e^1+i e^2.
\end{gather*}
Here $\left\{\left(\frac{\p}{\p u^1}\right)_u, \left(\frac{\p}{\p
u^2}\right)_u \right\}$ is the standard basis on $T_uN$. Then it
follows from \eqref{equation:frame1} that, $\zeta$ satisf\/ies
\begin{gather}
\zeta_x + 2 (\log \lambda )_z z_x \zeta=0. \label{equation:local1}
\end{gather}
The equation \eqref{equation:local1} is a linear f\/irst-order
ordinary dif\/ferential equation, and
\eqref{equation:frame1}--\eqref{equation:bound} has a global
solution $e$.

Since $\left\{e, Je   \right\}$ is an orthonormal frame along $u$,
we have
\begin{gather*}
0=\p_t [g_u(e,e)]=2g_u(\nabla_te,e),
\end{gather*}
which implies
\begin{gather}
\nabla_te=\alpha Je, \label{equation:frame2}
\end{gather}
for some real-valued function $\alpha (t,x)$. Moreover, in this
frame the coordinates of $u_t$ and $u_x$ are given by two complex
valued functions $p$ and $q$, namely, set
\begin{gather}
u_t =p_1e+p_2Je, \label{equation:frame3}
\\
u_x =q_1e+q_2Je, \label{equation:frame4}
\end{gather}
where $p_1$, $p_2$, $q_1$, $q_2$ are real-valued functions of
$(t,x)$, and set
\begin{gather*}
p=p_1+ip_2, \qquad \text{and} \qquad q=q_1+iq_2.
\end{gather*}
The relationships between $p$, $q$ and $\alpha$ are as follows.
Using the equation \eqref{equation:pde3}, we have
\begin{gather}
p=aq_{xx}+iq_x+b\lvert{q}\rvert^2q. \label{equation:change1}
\end{gather}
The compatibility condition $\nabla_xu_t=\nabla_tu_x$ implies
\begin{gather}
p_x=q_t+i\alpha q. \label{equation:change2}
\end{gather}
Combining \eqref{equation:change1} and \eqref{equation:change2},
we have
\begin{gather}
q_t=aq_{xxx}+iq_{xx}-i\alpha q +b\big(\left| q \right|^2q \big)_x.
\label{equation:change3}
\end{gather}
Moreover, we can obtain
\begin{gather}
\alpha_x=-\kappa(u)\operatorname{Im}\,(\bar{q}p),
\label{equation:change4}
\end{gather}
where $\kappa(u)=g_u(R(e,Je)Je,e)$ is the Gaussian curvature of
$N$ at $u$, and $R$ is the Riemann curvature tensor. Indeed, on
one hand, using \eqref{equation:frame1}, \eqref{equation:frame2}
and the identity
$\nabla_t\nabla_xe-\nabla_x\nabla_te=R(u_t,u_x)e$, we have
\begin{gather}
R(u_t,u_x)e=-\nabla_x(\alpha Je)=-\alpha_xJe. \nonumber
\end{gather}
On the other hand, combining \eqref{equation:frame3},
\eqref{equation:frame4} and the identity $R(X,Y)e=-R(Y,X)e$ for
any  vector f\/ields $X$ and $Y$ along $u$, we have
\begin{gather*}
R(u_t,u_x)e = (p_2q_1-p_1q_2)R(Je, e)e =
\operatorname{Im}\,(\bar{q}p)R(Je,e)e.
\end{gather*}
Hence it follows that
\begin{gather}
-\alpha_xJe=\operatorname{Im}\,(\bar{q}p)R(Je,e)e.
\label{equation:yanagiharakanako}
\end{gather}
Thus taking the inner product of \eqref{equation:yanagiharakanako}
and $-Je$, we have \eqref{equation:change4}.

Substituting \eqref{equation:change1} into
\eqref{equation:change4}, we have
\begin{gather*}
\alpha_x= \kappa(u) \left( -\frac{1}{2}\left|q\right|^2 + a
\operatorname{Im}\,(q\bar{q}_x) \right)_x.
\end{gather*}
Integrating this over $(-\infty,x]$, we get
\begin{gather}
\alpha(t,x)=
  \int_{-\infty}^x
  \kappa(u)(t,x')\left(
  -
  \frac{1}{2}\left|q\right|^2 + a \operatorname{Im}\,(q\bar{q}_x)
  \right)_x(t,x')dx'
  +
  \alpha(t,-\infty)
\nonumber
\\
\phantom{\alpha(t,x)}{} =
  \kappa(u)(t,x)\left(
  -
  \frac{1}{2}\left|q\right|^2 + a \operatorname{Im}\,(q\bar{q}_x)
  \right)(t,x)
+
  \frac{1}{2}
  \int_{-\infty}^x
  \left(\kappa(u) \right)_x(t,x')
  \left|q(t,x')\right|^2
  dx'
\nonumber
\\
\phantom{\alpha(t,x)=}{}  -
  a
  \int_{-\infty}^x
  \left(\kappa(u) \right)_x(t,x')
  \operatorname{Im}\,(q\bar{q}_x)(t,x')dx'
  +
  A(t),
\label{equation:nanako}
\end{gather}
where $A(t)=\alpha(t, -\infty)$. Substituting
\eqref{equation:nanako} into \eqref{equation:change3}, we have
\begin{gather*}
  q_t -aq_{xxx}-iq_{xx}
 =
  \left(\frac{a}{2}\kappa(u)+2b\right)
  \lvert{q}\rvert^{2}
  q_{x}
  -
  \left(\frac{a}{2}\kappa(u)-b\right)
  q^{2}\bar{q}_{x}
\\
\qquad{}
  +
  ia
  \left[
  \int_{-\infty}^{x}
  (\kappa(u))_{x}\operatorname{Im}\,(q\bar{q}_{x})
  dx'
  \right]
  q
  -
  \frac{i}{2}
  \left[
  \int_{-\infty}^{x}
  (\kappa(u))_{x}
  \lvert{q}\rvert^{2}
  dx'
  \right]
  q
  +
  \frac{i}{2}
  \kappa(u)\lvert{q}\rvert^{2}q
  -
  iA(t)q.
\end{gather*}
If we set
\begin{gather}
Q(t,x)=q(t,x)\exp\left(i\int_0^tA(\tau)d\tau \right),
\label{equation:ebichan}
\end{gather}
$Q(t,x)$ solves \eqref{equation:t3rd}. This completes the proof of
Theorem~\ref{theorem:hasimoto3}.
\end{proof}

\begin{proof}[Proof of Theorem~\ref{theorem:hasimoto4}]
In the same way as the previous proof, construct the moving frame
$\left\{ e, Je\right\}$ and set complex-valued functions $p$ and
$q$.

The relationships between $p$, $q$ and $\alpha$ are as follows.
The equation \eqref{equation:pde4} shows that
\begin{gather}
-ip=q_x-a q_{xxx}+ b|q|^2q_x+ \frac{c}{2}\left( |q|^2\right)_xq.
\label{equation:change21}
\end{gather}
Using $\nabla_xu_t=\nabla_tu_x$, we have
\begin{gather}
p_x=q_t+i\alpha q. \label{equation:change22}
\end{gather}
Combining \eqref{equation:change21} and \eqref{equation:change22},
we have
\begin{gather}
-iq_t=-a q_{xxxx}+q_{xx}+b\left( |q|^2q_x \right)_x +\frac{c}{2}
\left( \left( |q|^2\right)_xq \right)_x-\alpha q.
\label{equation:change23}
\end{gather}
Using the curvature tensor
$\nabla_t\nabla_xe-\nabla_x\nabla_te=R(u_t,u_x)e$, we have
\begin{gather}
\alpha_x=-\kappa(u)\operatorname{Im}\,(\bar{q}p),
\label{equation:change24}
\end{gather}
where $\kappa(u)=g_u(R(e,Je)Je,e)$ is the Gaussian curvature of
$N$ at $u$.

Substituting \eqref{equation:change1} into
\eqref{equation:change24}, we have
\begin{gather}
\alpha_x = - \kappa(u) \left( \frac{1}{2}\left|q\right|^2 +
\frac{b+c}{4}\left|q\right|^4 \right)_x + \frac{a}{2}\kappa(u)
(\bar{q}_{xxx}q+q_{xxx}\bar{q}). \label{equation:mao}
\end{gather}
Substituting the integration of the right hand side of
\eqref{equation:mao} on $(-\infty,x]$ into
\eqref{equation:change23}, we have
\begin{gather}
 q_t+ia q_{xxxx}-iq_{xx}
  =
 ib
  \bigl(\lvert{q}\rvert^2q_x\bigr)_x
  +
  i\frac{c}{2}
  \bigl\{
  \bigl(\lvert{q}\rvert^2\bigr)_x
  q
  \bigr\}_x
  -
  iA(t)q
\nonumber
\\
\qquad{} -
  iq
  \int_{-\infty}^x
  \left\{
  -
  \kappa(u)
  \left(
  \frac{1}{2}\lvert{q}\rvert^2
  +
  \frac{b+c}{4}\lvert{q}\rvert^4
  \right)_x
  +
  \frac{a}{2}\kappa(u)
  (\bar{q}_{xxx}q+q_{xxx}\bar{q})
  \right\}
  dx^\prime.
\label{equation:tadanohitoshi}
\end{gather}
Calculation of the fourth terms of the right hand side of
\eqref{equation:tadanohitoshi} by the use of integrations by parts
implies that
\begin{gather*}
 q_t +ia q_{xxxx}-iq_{xx}
  =
 i
  \left(
  b
  +
  \frac{c}{2}
  -
  \frac{a}{2}\kappa(u)
  \right)
  \lvert{q}\rvert^2
  q_{xx}
  +
  i
  \left(
  \frac{c}{2}
  -
  \frac{a}{2}\kappa(u)
  \right)
  q^2\bar{q}_{xx}
  +
  i
  \left(
  b
  +
  \frac{c}{2}
  \right)
  \bar{q}q_{x}^2
\\
\qquad{} +i
  \left(
  b
  +
  \frac{3c}{2}
  +
  \frac{a}{2}\kappa(u)
  \right)
  q\lvert{q_x}\rvert^2
  +
  \frac{i}{2}
  \kappa(u)
  \lvert{q}\rvert^2q
  +i
  \frac{b+c}{4}
  \kappa(u)
  \lvert{q}\rvert^4q
\\
\qquad{} -
  i
  \left[
  \int_{-\infty}^x
  \Bigl(\kappa(u)\Bigr)_x
  \left\{
  \frac{1}{2}
  \lvert{q}\rvert^2
  +
   \frac{b+c}{4}
  \lvert{q}\rvert^4
  +
  \frac{a}{2}
  \lvert{q_x}\rvert^2
  -
  \frac{a}{2}
  q_{xx}\bar{q}
  -
   \frac{a}{2}
  \bar{q}_{xx}q
  \right\}
  dx'
  \right]
  q
-iA(t)q.
\end{gather*}
We can take $A(t)\equiv 0$ without loss of generality by the same
way as we use \eqref{equation:ebichan} in the previous proof. That
is, if we set
\begin{gather*}
Q(t,x)=q(t,x)\exp\left(i\int_0^tA(\tau)d\tau \right),
\end{gather*}
it is easy to check that $Q$ solves \eqref{equation:t4th}. Thus we
complete the proof.
\end{proof}

\section{Concluding remarks}
\label{section:remarks}

Our reduction via moving frame leads to an understanding of the
essential structure of the dif\/ferential equations
\eqref{equation:pde3} and \eqref{equation:pde4}. Unfortunately,
however, the results of this paper alone are not~suf\/f\/icient
for establishing well-posedness of the initial-value problem. In
this section we discuss the relationship between our reduction and
the theory of linear partial dif\/ferential ope\-rators. One can
consult with \cite{Mizohata} on the basic theory of linear
dispersive partial dif\/ferential operators.

Here we review the theory of one-dimensional linear dispersive
partial dif\/ferential operators. Consider the initial value
problem for one-dimensional Schr\"odinger-type evolution equations
of the form
\begin{gather}
  v_t-iv_{xx}+b(x)v_x+c(x)v
 =
  f(t,x)
  \qquad\text{in}\quad
 \mathbb{R}\times\mathbb{R},
\label{equation:schroedinger}
\\
  v(0,x)
 =
  v_0(x)
  \qquad\text{in}\quad
 \mathbb{R},
\label{equation:data}
\end{gather}
where $b(x),c(x)\in\mathscr{B}^{\infty}(\mathbb{R})$, which is the
set of all smooth functions on $\mathbb{R}$ whose derivative of
any order are bounded on $\mathbb{R}$, $v(t,x)$ is a
complex-valued unknown function, and $f(t,x)$ and $v_0(x)$ are
given functions. It is well-known that the initial value problem
\eqref{equation:schroedinger}--\eqref{equation:data} is well-posed
in the framework of $L^2(\mathbb{R})$ if and only if
\begin{gather}
\sup_{x\in\mathbb{R}} \left\lvert \int_0^x
\operatorname{Im}\,b(y)dy \right\rvert < \infty.
\label{equation:mizohata}
\end{gather}
See \cite{Mizohata} for details. This condition describes the
balance of the local smoothing ef\/fect of $\exp(it\p_x^2)$
\begin{gather}
\lVert (1+x^2)^{-\delta} (-\p_x^2)^{1/4} e^{it\p_x^2}\phi
\rVert_{L^2(\mathbb{R}^2)} \leqslant C
\lVert\phi\rVert_{L^2(\mathbb{R})}, \qquad \delta>1/4, \nonumber
\end{gather}
and the strength of the bad part of the f\/irst order term
$i\operatorname{Im}\,\{b(x)\}\p_x$. See, e.g., \cite{CS,Sjolin}
and \cite{Vega} for the local smoothing ef\/fect of linear
dispersive partial dif\/ferential equations with constant
coef\/f\/icients. If the condition~\eqref{equation:mizohata}
breaks down, the equation \eqref{equation:schroedinger} behaves
like the Cauchy--Riemann equation $v_i+iv_x=0$ and the initial
value problem \eqref{equation:schroedinger}--\eqref{equation:data}
can be solvable only in some framework of real-analytic functions.
Under the condition \eqref{equation:mizohata}, a gauge transform
\begin{gather}
v(t,x) \longmapsto w(t,x) = v(t,x) \exp \left( \frac{i}{2}
\int_0^x b(y)dy \right) \label{equation:gauge}
\end{gather}
is automorphic on $C(\mathbb{R};L^2(\mathbb{R}))$, and
\eqref{equation:schroedinger}--\eqref{equation:data} becomes the
initial value problem of the form
\begin{gather*}
w_t - iw_{xx} + d(x)w = g(t,x), \qquad w(0,x)=w_0(x),
\end{gather*}
which is easy to solve. We remark that if we replace
$x\in\mathbb{R}$ by $x{\in}\mathbb{R}/\mathbb{Z}$, then the local
smoothing ef\/fect breaks down and the condition
\eqref{equation:mizohata} becomes
\begin{gather}
\int_0^1\operatorname{Im}\,b(y)dy=0. \label{equation:mizohataT}
\end{gather}
The generalized Hasimoto transform introduced in \cite{CSU}
corresponds to \eqref{equation:gauge}. The equation
\eqref{equation:pde3} with $a=b=0$ behaves like symmetric
hyperbolic systems since $\nabla^NJ=0$, where $\nabla^N$ is the
Levi-Civita connection of $(N,J,g)$. This fact is implicitly
applied to proving of the existence of solutions to the initial
value problem. See \cite{Koiso,Koiso2,PWW} and~\cite{SSB}. From a
point of view of the theory of linear partial dif\/ferential
operators, this is due to the fact that the $L^2$-norm given by
the K\"ahler metric of the form
\begin{gather}
\left\{ \int_\mathbb{R} g(V,V) dx \right\}^{1/2}
\qquad\text{for}\quad V{\in}\Gamma_0(u^{-1}TN),
\label{equation:chiemi}
\end{gather}
corresponds to the modif\/ied norm
\begin{gather*}
\left\{ \int_\mathbb{R} \lvert{v}\rvert^2 \exp \left(- \int_0^x
\operatorname{Im}\,b(y)dy \right) dx \right\}^{1/2}
\end{gather*}
for \eqref{equation:schroedinger}. Here $\Gamma_0(u^{-1}TN)$ is
the set of compactly supported smooth sections of $u^{-1}TN$.

Similarly, we discuss one-dimensional third and fourth-order
linear dispersive partial dif\/fe\-ren\-tial equations of the forms
\begin{gather}
   v_t+v_{xxx}+\alpha(x)v_{xx}+\beta(x)v_x+\gamma(x)v
 =0,
\label{equation:tarama}
\\
  v_t+v_{xxxx}+A(x)v_{xxx}+B(x)v_{xx}+C(x)v_x+D(x)v
 =0,
\label{equation:mizuhara}
\end{gather}
where $v:\mathbb{R}\times\mathbb{R}\rightarrow\mathbb{C}$ is an
unknown function, and
$\alpha,\beta,\gamma,A,B,C,D\in\mathscr{B}^\infty(\mathbb{R})$.
Tarama studied the initial value problem for
\eqref{equation:tarama} in \cite{Tarama1} and \cite{Tarama2}, and
Mizuhara studied the initial value problem for
\eqref{equation:mizuhara} in \cite{Mizuhara}. The initial value
problem for \eqref{equation:tarama} is $L^2$-well-posed if and
only if
\begin{gather}
\sup_{x\in\mathbb{R}} \left\lvert \int_0^x
\operatorname{Re}\,\alpha(s)ds \right\rvert <\infty,
\label{equation:oogaki}
\\
\sup_{\substack{(x,y)\in\mathbb{R}^2 \\ x\neq y}} \left\lvert
\int_y^x \operatorname{Im}\, \left\{
\beta(s)+\frac{\alpha(s)^2}{3} \right\} ds \right\rvert
\lvert{x-y}\rvert^{-1/2} < \infty. \label{equation:yuko}
\end{gather}
Under the conditions \eqref{equation:oogaki} and
\eqref{equation:yuko}, the second order term in
\eqref{equation:tarama} is eliminated by a gauge transform
\begin{gather}
v(t,x) \longmapsto w(t,x) = v(t,x) \exp \left( \frac{1}{3}
\int_0^x \alpha(y)dy \right), \label{equation:kojimayoshio}
\end{gather}
and the bad part of the f\/irst order term of
\eqref{equation:tarama} is canceled out by a transform def\/ined
by a~pseudodif\/ferential operator. See \cite{Tarama1} and
\cite{Tarama2} for details. When $x\in\mathbb{R}/\mathbb{Z}$,
\eqref{equation:oogaki} and \eqref{equation:yuko} become
\begin{gather*}
\int_0^1 \operatorname{Re}\,\alpha(y)dy = \int_0^1
\operatorname{Im}\, \left\{ \beta(y)+\frac{\alpha(y)^2}{3}
\right\}dy =0.
\end{gather*}
The generalized Hasimoto transform in the proof of
Theorem~\ref{theorem:hasimoto3} corresponds to
\eqref{equation:kojimayoshio}. In~\cite{Onodera} the author proved
local and global existence theorems of the initial value problem
for \eqref{equation:pde3} on
$\mathbb{R}\times\mathbb{R}/\mathbb{Z}$ with values in a K\"ahler
manifold. He made use of the classical energy estimates of
\eqref{equation:chiemi}. Notice that \eqref{equation:pde3} also
behaves like symmetric hyperbolic systems since $\nabla^Ng=0$ and
$\nabla^NJ=0$. The $L^2$-norm \eqref{equation:chiemi} works as if
the unknown functions were transformed by
\eqref{equation:kojimayoshio} and the pseudodif\/ferential
operator. In other words, if one see the metric tensor $g$ as a
real symmetric matrix or a hermitian matrix, a commutator
$[g,I\p_{xxx}+J\p_{xx}]$ eliminates the bad part of the second and
the f\/irst order terms, where $I$ is the identity matrix. It is
worth to mention that Nishiyama and Tani studied existence
theorems of the initial value problem for the third-order
two-sphere valued model \eqref{equation:vortex9} in \cite{NT} and
\cite{TN}.

The necessary and suf\/f\/icient condition of the
$L^2$-well-posedness of the initial value problem for
\eqref{equation:mizuhara} was given in \cite{Mizuhara} under some
technical condition. One cannot expect the existence of a~smooth
closed curve, that is a periodic solution in $x$, since the
classical energy method breaks down for \eqref{equation:pde4}.

\subsection*{Acknowledgements}
The author expresses gratitude to Hiroyuki Chihara for several
discussions and valuable advice. Also, thanks to the referees for
carefully reading the manuscript. The author is supported by the
JSPS Research Fellowships for Young Scientists and the JSPS
Grant-in-Aid for Scientif\/ic Research No.19$\cdot$3304.

\pdfbookmark[1]{References}{ref}
\LastPageEnding


\begin{thebibliography}{99}

\footnotesize\itemsep=0pt

\bibitem{CSU}
Chang~N.H., Shatah~J., Uhlenbeck~K.,
 Schr\"odinger maps,
{\it Comm.  Pure Appl. Math.} {\bf 53} (2000), 590--602.

\bibitem{CS}
Constantin~P., Saut~J.C.,
 Local smoothing properties of dispersive equations,
{\it J. Amer. Math. Soc.} {\bf 1} (1989), 413--446.

\bibitem{DR}
Da~Rios,
 On the motion of an unbounded f\/luid
with a vortex f\/ilament of any shape, {\it Rend.  Circ. Mat.
Palermo} {\bf 22} (1906), 117--135 (in Italian).

\bibitem{Fukumoto}
Fukumoto Y.,
 Motion of a curved vortex f\/ilament: higher-order asymptotics,
in Proc. of IUTAM Symposium on Geometry and Statistics of
Turbulence (1999, Hayama), Editors
 T.~Kambe, T.~Nakano and T.~Miyauchi, {\it Fluid Mech. Appl.}, Vol.~59, Kluwer Acad. Publ., Dordrecht, 2001, 211--216.


\bibitem{FM}
Fukumoto~Y., Miyazaki~T.,
 Three-dimensional distortions of a vortex f\/ilament with axial velocity,
{\it J. Fluid Mech.} {\bf 222} (1991), 369--416.

\bibitem{Hama}
Hama~F.R.,
 Progressive deformation of a curved vortex f\/ilament
by its own induction, {\it Phys. Fluids} {\bf 5} (1962),
1156--1162.

\bibitem{Hasimoto}
Hasimoto H.,
 A soliton on a vortex f\/ilament,
{\it J. Fluid. Mech.} {\bf 51} (1972), 477--485.

\bibitem{Koiso}
Koiso N.,
 The vortex f\/ilament equation and
a semilinear Schr\"odinger equation in a Hermitian symmetric
space, {\it Osaka J. Math.} {\bf 34} (1997), 199--214.

\bibitem{Koiso2}
Koiso N.,
 Vortex f\/ilament equation and
 semilinear Schr\"odinger equation,
 in Nonlinear Waves (1995, Sapporo),
{\it Gakuto Intern. Ser. Math. Sci. Appl.}, Vol.~10, 1997,
231--236.

\bibitem{Mizohata}
Mizohata M., On the Cauchy problem, Academic Press, 1985.

\bibitem{Mizuhara}
Mizuhara R.,
 The initial value problem for
third and fourth order dispersive equations in one space
dimensions, {\it Funkcial. Ekvac.} {\bf 49} (2006), 1--38.

\bibitem{NT}
Nishiyama T., Tani A.,
 Initial and initial-boundary value problems for a vortex f\/ilament
    with or without axial f\/low,
{\it SIAM J. Math. Anal. } {\bf 27} (1996), 1015--1023.

\bibitem{Onodera}
Onodera E.,
 A third-order dispersive f\/low for closed curves into
K\"ahler manifolds, {\it J. Geom. Anal.} {\bf 18} (2008), to
appear, \href{http://arxiv.org/abs/0707.2660}{arXiv:0707.2660}.

\bibitem{PWW}
Pang P.Y.Y., Wang H.Y., Wang~Y.D.,
  Schr\"odinger f\/low on Hermitian locally symmetric spaces,
 {\it Comm. Anal. Geom.} {\bf 10} (2002), 653--681.

\bibitem{Sjolin}
Sj\"olin P.,
 Regularity of solutions to the Schr\"odinger equation,
{\it Duke  Math. J.}  {\bf 55} (1987), 699--715.


\bibitem{SSB}
Sulem P.-L., Sulem C., Bardos C.,
 On the continuous limit for a system of classical spins,
{\it Comm. Math. Phys.} {\bf 107} (1986), 431--454.

\bibitem{TN}
Tani A., Nishiyama T.,
 Solvability of equations for motion of a vortex f\/ilament
with or without axial f\/low, {\it Publ. Res. Inst. Math. Sci.}
{\bf 33} (1997), 509--526.

\bibitem{Tarama1}
Tarama S.,
 On the wellposed Cauchy problem for some dispersive equations,
{\it J. Math. Soc. Japan} {\bf 47} (1995), 143--158.

\bibitem{Tarama2}
Tarama S.,
 Remarks on $L\sp2$-wellposed Cauchy problem for some
dispersive equations, {\it J. Math. Kyoto Univ.} {\bf 37} (1997),
757--765.

\bibitem{Vega}
Vega L.,
 The Schro\"odinger equation: pointwise convergence to the initial date,
{\it Proc. Amer.  Math. Soc.} {\bf 102} (1988), 874--878.

\end{thebibliography}
\end{document}